\documentclass{svmult}

\usepackage{makeidx}
\usepackage{graphicx}
\usepackage{multicol}
\usepackage{amsfonts}

\newcommand{\R}{{\mathbb{R}}}
\newcommand{\Lipp}{{\mathrm{Lip}'}}

\makeindex

\begin{document}

\title*{
   Convergence of approximate solutions of conservation laws
   }
\author{
   Sebastian Noelle\inst{1}\and
   Michael Westdickenberg\inst{2}
   \thanks{This work was supported by SFB 256 at Bonn University.}}
\institute{
   Institut f\"ur Geometrie und Praktische Mathematik \\
   Rheinisch--Westf\"alische--Technische Hochschule Aachen \\
   Templergraben 55, 52056 Aachen, Germany \\
      \textit{noelle@igpm.rwth-aachen.de}
\and
   Institut f\"ur Angewandte Mathematik \\
   Rheinische Friedrich-Wilhelms-Universit\"at Bonn \\
   Wegelerstrasse 10, 53115 Bonn, Germany \\
      \textit{mwest@iam.uni-bonn.de}}

\titlerunning{
   Convergence of approximate solutions of conservation laws}

\maketitle

\begin{abstract}
In this paper we consider convergence of approximate solutions
of conservation laws. We start with an overview over the
historical developments since the 1950s, and the analytical tools
used in this context. Then we present some of our own results
on the convergence of numerical approximations, discuss recent
related work and open problems.
\end{abstract}

\section{Exact and approximate solutions of conservation laws}
\label{Intro}
Hyperbolic conservation laws are used widely to model the local
conservation of physical quantities like density, momentum and energy.
Applications include gasdynamics, ideal magnetohydrodynamics,
shallow water and traffic flows. The equations take the form
\begin{equation}
\label{1_100}
\frac{\partial u}{\partial t}+\nabla\cdot f(u)=0,
\quad\mbox{for $(t,x)\in[0,T]\times\R^d$},
\end{equation}
where $u(t,x):=(u_1(t,x),\ldots,u_m(t,x))$ is the $m$-vector of
conservative variables with initial data $u(0,\cdot)=:u_0$, and
$f(u)=(f_1(u),\ldots,f_d(u))$ is the flux function.
It is well known that for nonlinear fluxes, even for smooth initial
data, the solution of (\ref{1_100}) may cease to exist in the
classical sense due to the formation of shocks in finite time
(cf. Riemann 1859 \cite{Riemann}).
Therefore it is necessary to consider weak solutions,
i.e. functions $u$ that are bounded and satisfy (\ref{1_100})
in the distributional sense.
Existence of weak solutions for general systems of conservation
laws is largely unknown, and weak solutions are not unique.
However, for scalar equations (with $m=1$) the Cauchy problem
is well understood.

As usual, a possible strategy to establish existence of weak
solutions of the conservation law (\ref{1_100}) is first to
regularize the problem and obtain a sequence of approximate
solutions $\{u^\epsilon\}_\epsilon$ for $\epsilon>0$, then to
show that a subsequence of $u^\epsilon$ converges to a limit
function $u$ which is a weak solution of (\ref{1_100}). The
regularization mostly used is the vanishing viscosity method,
where $u^\epsilon$ solves the parabolic problem $\partial_t 
u^\epsilon+\nabla\cdot f(u^\epsilon)=\epsilon\Delta u^\epsilon$.
Since in general exact solutions of conservation laws 
are not known, approximations through numerical computations
are of utmost importance
in the applications. As the grid size tends to zero, again
a sequence of approximate solutions is created,
and then a convergence analysis is important to make sure
that the results are reliable. We say more on that issue
in Section \ref{S3}.

We will assume that the sequence $\{u^\epsilon\}_\epsilon$ of 
approximate solutions is uniformly bounded in $L^\infty([0,T]
\times\R^d)$. For the vanishing viscosity method and for many 
numerical schemes this can be obtained from a (discrete) maximum 
principle or an invariant region argument. Then it is 
possible to extract a subsequence which converges weak* to
some limit function $u\in L^\infty([0,T]\times\R^d)$. But in
order to prove that $u$ is a weak solution of (\ref{1_100}),
it is necessary to show that the weak limit of $f(u^\epsilon)$
coincides with $f(u)$. For a weakly converging sequence
$u^\epsilon$ and a nonlinear flux $f$ this is not true in general because
oscillations may occur. A sufficient condition, however, would
be strong $L^1_{loc}$-convergence. So the question is how this
additional information can be obtained.

We already mentioned that for a general nonlinear hyperbolic
conservation law there may exist many different weak solutions
corresponding to the same initial data. In order to ensure
uniqueness one imposes an additional condition, often called
an entropy condition because of its analogy with the second law
of thermodynamics, which selects the physically relevant weak
solution out of all possible ones. For scalar conservation laws
several related formulations are in use. We give the one that
relies on convex entropies.

We call a pair of functions $\eta\colon\R^m\rightarrow\R$ and
$q\colon\R^m\rightarrow\R^d$ an entropy-entropy flux pair if
$\eta$ is convex and if the compatibility relations
\begin{equation}
\label{1_150}
q_k'(u)=\eta'(u)f_k'(u)
\quad\mbox{for all $u\in\R^m$ and $k=1,\ldots,d$}
\end{equation}
hold.
Then the entropy condition reads as follows: Out of all weak
solutions of the conservation law (\ref{1_100}) given above,
we only select those satisfying
\begin{equation}
\label{1_200}
\frac{\partial\eta(u)}{\partial t}+\nabla\cdot q(u)\leq 0
\quad\mbox{in distributional sense}
\end{equation}
for all entropy-entropy flux pairs $(\eta,q)$. We also assume that
the initial entropy is bounded: $\int_{\R^d}\eta(u_0)\,dx<\infty$.
These solutions are called weak entropy solutions.
For scalar conservation laws all convex functions $\eta$ can 
serve as entropies, and it turns out that weak entropy solutions
are unique.
This is exceptional. For systems of conservation laws there often
exists only a limited number or even only one single (the physical)
entropy, because the corresponding compatibility relations are much
more restrictive.
The entropy inequalities provide crucial estimates in the
convergence analyis for $\{u^\epsilon\}_\epsilon$, and the
lack of enough entropy-entropy flux pairs is the main reason
why the existence theory for systems of conservation laws
is not as well developed as for the scalar case.
For more on the entropy condition
we refer to Lax \cite{Lax1971} and Liu \cite{Liu}.
Further information on hyperbolic conservation laws can 
be found in Lax \cite{Lax1973}, Smoller \cite{Smoller} and
Dafermos \cite{Dafermos}.

Our paper is organized as follows.
In Section 2, we give an overview over the developments in the
theory of conservation laws 
since the 1950s, and the analytical tools used to study
the convergence of approximate solutions.
These include classical compactness arguments, approximation theory,
weak convergence methods (compensated compactness and measure-valued
solutions), and the kinetic approach. In Section 3, we present some
of our own work on the convergence of numerical approximations.
Finally, in Section 4, we discuss recent related work and open problems.

\section{Historical remarks on compactness arguments}
In this section we give an overview over some analytical tools
used to study convergence of sequences of approximate
solutions for conservation laws.
\subsection{Regularity estimates and compactness}
\label{S21}
The classical approach to prove strong $L^1_{loc}$-convergence of a
sequence $\{u^\epsilon\}_\epsilon$ of functions is to prove that the
$u^\epsilon$ possess some positive regularity, i.e. the sequence is
uniformly bounded in some appropriate function space. If the function
space is $BV$, the space of functions with bounded variation, then 
strong compactness in $L^1_{loc}$ follows from Helly's Theorem. 

In 1957 Ladyzhenskaya \cite{Ladyzhenskaya} and Oleinik \cite{Oleinik1957},
both former students of Petrovskii in Moscow, independently published
results on a scalar conservation law with strictly convex 
flux in one space dimension. The entropy condition they used to select
the physical solution (today known as ``condition $E$'') can be
interpreted as a one-sided Lipschitz 
condition: For any time $t>0$, the solution $u(t,\cdot)$
is Lipschitz continuous where it is increasing (in $x$). There
is no restriction where it is decreasing. In particular, decreasing jump 
discontinuities are allowed. From this a $BV$ bound can be derived.
We would like to remark that the one-sided Lipschitz continuity was
already observed by E. Hopf \cite{Hopf} in his pioneering paper on
the Burgers equation, but he did not use it systematically as 
selection criterion.
While Ladyzhenskaya and Oleinik both studied the vanishing viscosity limit,
Oleinik together with her student Vvedenskaya \cite{Vvedenskaya1956,Oleinik1957}
also proved convergence of the Lax-Friedrichs scheme \cite{Lax1954}. 
In 1963, Oleinik's paper was translated into English, while Ladyzhenskaya
and Vvedenskaya's work was not. Oleinik's contribution is also discussed in full
detail in Smoller's famous book \cite{Smoller}. This may
explain why condition $E$ is nowadays usually attributed to
Oleinik.

The multidimensional scalar problem was studied by Conway and Smoller 
\cite{ConwaySmoller} and Vol'pert \cite{Volpert} in the 1960s, and by Kruzkov 
\cite{Kruzkov} in 1970 who proved existence and uniqueness
of weak entropy solutions for the scalar conservation law. Kruzkov's
main observation was that the solution operator of the scalar conservation
law is an $L^1$-contraction: He used the special convex entropies
\begin{equation}
\eta(u,k):=|u-k|
\quad\mbox{and}\quad
q(u,k):=\mathrm{sgn}(u-k)|f(u)-f(k)|
\end{equation}
(nowadays known as Kruzkov entropies), where $k\in\R$ is a real parameter, 
and a ``doubling of variables''-argument to conclude that 
\begin{equation}
\label{21_2400}
\int_{\R^d}|u(t,x)-v(t,x)|dx
\leq\int_{\R^d}|u_0(x)-v_0(x)|dx
\end{equation}
for all times $t>0$. Here  $u$ and $v$ are weak entropy solutions 
of (\ref{1_100}) corresponding to initial data $u_0$ and $v_0$,
respectively. Uniqueness of weak entropy solutions of (\ref{1_100}) is 
an immediate consequence of this estimate. And since the solution
operator of (\ref{1_100}) is invariant under translations in
$x$, we can apply (\ref{21_2400}) also to $v(t,x):=u(t,x+h)$
for any $h\in\R^d$ and obtain from this a $BV$-bound $\|u(t,
\cdot)\|_{BV(\R^d)}\leq\|u_0\|_{BV(\R^d)}$.
An analogous estimate for temporal regularity can then be 
derived using (\ref{1_100}).
So Helly's theorem can be applied.

Let us briefly mention that there is also a well-developed theory 
for systems of conservation laws in one space dimension relying on
$BV$ bounds. This theory started with the work of Glimm \cite{Glimm}
in 1965 who introduced a numerical scheme that produces approximate 
solutions with bounded variation. Later on, Bressan and his coworkers 
extended Glimm's ideas considerably and set up a very powerful 
framework for studying existence and uniqueness of systems of 
conservation laws in 1-d. We refer to Bressan \cite{Bressan}.

$BV$ bounds also play a major role in the design of higher
order accurate numerical algorithms. In 1983, Harten \cite{Harten}
introduced the class of TVD (total variation diminishing) schemes
and proved their convergence in the linear case. Osher \cite{Osher}
introduced the class of ``$E$-schemes'', the most general class
of schemes which satisfy a discrete entropy inequality,
and proved convergence of second order accurate
semidiscrete TVD schemes for nonlinear scalar conservation laws.
Convergence of a second order accurate fully discrete scheme was
proved by Osher and Tadmor in 1988 \cite{OsherTadmor}.

\subsection{Approximation theory}

Approximation theory in the context of scalar conservation laws
means to estimate the distance between an approximate solution 
(at some time $t>0$, say) and the uniquely defined weak entropy
solution which is already known to exist. This approach gives
error estimates and convergence rates. Applied to numerical
schemes it tells you how close an approximate solution is to
the exact one. 

The first who gave error estimates for approximate
solutions of scalar conservation laws in several space dimensions
was Kuznetsov \cite{Kuznetsov}, who studied the vanishing viscosity 
method and the Lax-Friedrichs scheme. His analysis relies on the 
$L^1$-contraction property of the solution operator of (\ref{1_100})
we discussed in the previous section. In fact, it is possible 
to obtain an estimate like (\ref{21_2400}) even if $v$ is not
an exact, but only a (suitable) approximate solution. Then 
additional terms measuring the approximation error enter on the 
right hand side which can be controlled if the initial data is
in $BV$. More precisely, one can estimate the
$L^1$-distance between the weak entropy solution $u$ and the
approximate solution $u^\epsilon$ at a given time $t$ in terms
of the distance and $BV$-norm of the initial data, the residual
etc. These estimates depend on $\epsilon$ (which
could be the gridsize in numerical schemes, for example),
and the convergence rate proved by Kuznetsov is
$\epsilon^{1/2}$.

Kuznetsov's ideas were used later on to give convergence rate 
estimates for more sophisticated numerical schemes, as well. We refer 
to Sanders \cite{Sanders}, Vila \cite{Vila}, Cockburn, Coquel and
LeFloch \cite{CockburnCoquelLeFloch}, Noelle \cite{Noelle1996}, and
Cockburn, Gremaud and Yang \cite{CockburnGremaudYang1998}.
We also refer to Bouchut and Perthame
\cite{BouchutPerthame} who reformulated Kruzkov and Kuznetsov's
approximation theory in a general, versatile form.

We also mention that recently, Tadmor \cite{Tadmor1991}
and coworkers developed a somewhat different approach
to proving error estimates and convergence rates for scalar
conservation laws in 1-d with strictly convex flux, 
the so-called $\Lipp$-theory. The idea is to measure the 
distance between the exact and the approximate solution not in 
the $L^1$-norm, but in the much weaker topology of $\Lipp(\R)$. 
This is the topological dual of the space of Lip\-schitz continuous 
functions. The approach can also give pointwise error estimates: The theory
shows that the convergence rate at some given point $x\in\R$ only 
depends on the regularity of the exact solution $u$ in a small 
neighborhood around $x$. 

\subsection{Weak convergence methods}
For modern numerical schemes using unstructured grids
it may be too hard or even impossible to prove
positive regularity in terms of boundedness of a sequence
$\{u^\epsilon\}_\epsilon$ in appropriate function spaces.
So the classical approach of 
proving compactness in $L^1_{loc}$ by using Helly's theorem or
Sobolev embeddings breaks down. There are other tools, which run
under the name of Weak Convergence Methods, that can prove strong
compactness without showing regularity first. These methods
are specialized in the sense that they rely heavily on the 
structure of the problem at hand. 

In 1954, Lax \cite{Lax1954} proved compactness in $L^1$ of the solution operator
of a one-dimensional, not necessarily strictly convex, scalar conservation law
using a weak topology. More than twenty years later, Tartar \cite{Tartar}
and Murat \cite{Murat} introduced an even more powerful method,
applicable also for systems of conservation laws,
the compensated compactness theory. A typical result is the
$\mathrm{div}$-$\mathrm{curl}$-Lemma: One assumes that
two sequences of vector-valued functions $\{U^\epsilon\}_\epsilon$ and 
$\{V^\epsilon\}_\epsilon$ are bounded in $(L^2_{loc}(\R^n))^N$ and that 
the divergence of $U^\epsilon$ and the curl of $V^\epsilon$ are both
precompact in $H^{-1}_{loc}(\R^n)$. Then there exist subsequences such that 
\begin{equation}
\label{22_2100}
U^\epsilon\cdot V^\epsilon\longrightarrow U\cdot V
\quad\mbox{in distributional sense}.
\end{equation}
Note that the assumptions on the differentiability of $U^\epsilon$ and
$V^\epsilon$ are not sufficient to obtain strong compactness in $L^2
_{loc}(\R^n)$ (from which (\ref{22_2100}) would follow trivially).
But in combination they give enough information to have weak continuity
of the scalar product.

The $\mathrm{div}$-$\mathrm{curl}$-Lemma has been used by Tartar to study the
convergence of sequences of approximate solutions $\{u^\epsilon\}_\epsilon$
of the scalar conservation law (\ref{1_100}) in one space dimension.
Compensated Compactness arguments can also be used to study certain 
systems of conservation laws in 1-d. As an example, we refer to the work
of DiPerna \cite{DiPerna1983} and Chen \cite{Chen} on isentropic Euler 
equations. 

Another weak convergence method is DiPerna's theory of measure-valued 
solutions for multidimensional scalar conservation laws.
As mentioned before, we assume that the sequence $\{u^\epsilon\}
_\epsilon$ of approximate solutions is uniformly bounded in $L^\infty$
and extract a subsequence converging weak* in $L^\infty$. Then the
sequence converges strongly in $L^1_{loc}$ if and only if no 
oscillations occur, i.e. if there is convergence pointwise a.e. 
The concept of measure-valued solution was first introduced by 
Tartar \cite{Tartar} in 1975, who used Young measures to 
describe the oscillations occuring in $\{u^\epsilon\}_\epsilon$.
A Young measure associated to $\{u^\epsilon\}_\epsilon$ is a 
weakly measurable mapping 
$\nu$ from $[0,T]\times\R^d$ into the space of probability measures 
$\mathrm{Prob}(\R)$, such that for all continuous functions 
$f$ we have $f(u^\epsilon)\rightarrow \, <\!\!\nu,f\!\!>$ in distributional sense, 
where the pairing $<\!\!\nu,f\!\!>$ is given by 
\[
<\!\nu,f\!>(t,x):=\int_\R f(\lambda) \,d\nu_{(t,x)}(\lambda).
\]
So the distributional limit of any nonlinear function of $u^\epsilon$
can be written down using one single $\nu$. 
A Young measure is called
a measure-valued solution of the Cauchy problem (\ref{1_100}) if
the following identity holds
\begin{equation}
\label{22_1000}
\frac{\partial}{\partial t}<\!\nu,\mathrm{id}\!>+\nabla\cdot<\!\nu,f\!>=0
\quad\mbox{in distributional sense}.
\end{equation}
Consistency with the entropy condition then means: For all convex
$\eta$
\begin{equation}
\label{22_1100}
\frac{\partial}{\partial t}<\!\nu,\eta\!>+\nabla\cdot<\!\nu,q\!>\leq 0
\quad\mbox{in distributional sense},
\end{equation}
where $q$ is the corresponding entropy flux.

Note that a measure-valued solution is even weaker than a weak
(distributional) solution. One can show that for any bounded sequence
of approximate solutions $\{u^\epsilon\}_\epsilon$ there exists a
Young measure $\nu$ and a subsequence converging to it in an
appropriate sense. That $\nu$ is a measure-valued solution
of (\ref{1_100}) then follows from the consistency of the approximation.

One has $<\!\!\nu,\mathrm{id}\!\!>=u$, where $u$ is the distributional
limit of $\{u^\epsilon\}_\epsilon$. But in general $<\!\!\nu,f\!\!>
\neq f(u)$. Equality holds for arbitrary nonlinear $f$ if and
only if
the Young measure reduces to a Dirac measure, i.e. if $\nu$ is of the form
$\nu_{(t,x)}=\delta_{u(t,x)}$ with $\delta_u$ the Dirac distribution
on $\R$ centered at $u$. This is equivalent to saying that the
sequence $\{u^\epsilon\}_\epsilon$ converges strongly in $L^1_{loc}$.

The following result is due to DiPerna \cite{DiPerna1985}
(see also Szepessy \cite{Szepessy1989a})
\begin{theorem}
\label{22_T100}
Assume that initial data $u_0\in L^1\cap L^\infty(\R^d)$ is given and
that there is exists a Young measure $\nu$ satisfying the following
properties \\
(i) the function $(t,x)\mapsto<\!\!\nu_{(t,x)},|\mathrm{id}|\!\!>$ is in
$L^\infty([0,T],L^1(\R^d))$, \\
(ii) $\nu$ is a measure-valued solution of (\ref{1_100}), \\
(iii) $\nu$ is consistent with the entropy condition, \\
(iv) $\nu$ assumes the initial data $u_0$ in the following sense
\begin{equation}
\label{22_1200}
\lim_{t\rightarrow 0}\frac 1t\int_0^t\int_{\R^d}<\!\nu_{(s,x)},
|\mathrm{id}-u_0(x)|\!>ds\,dx = 0.
\end{equation}
Then the Young measure reduces to a Dirac measure, i.e. $\nu_{(t,x)}
=\delta_{u(t,x)}$ for a.e. $(t,x)\in[0,T]\times\R^d$, where $u$ is
the unique entropy solution of (\ref{1_100}).
\end{theorem}

Assumption (iv) says that the initial data must be attained
in a stronger sense than just the sense of distributions.
It excludes the occurence of oscillations
in the sequence of approximate initial data $\{u_0^\epsilon
\}_\epsilon$. Then condition (iii) assures that no oscillations
can develop at later times: the only oscillations that can exist in the
sequence $\{u^\epsilon\}_\epsilon$ are those transported into the system
from the initial data. If the initial data converge strongly, so does
$\{u^\epsilon\}_\epsilon$ at any later time.

Theorem \ref{22_T100} was used to prove convergence of the streamline
diffusion shock-capturing method by Szepessy \cite{Szepessy1989b},
finite difference methods by Coquel and LeFloch \cite{CoquelLeFloch},
spectral viscosity approximations by Chen, Du and Tadmor \cite{ChenDuTadmor}
and finite volume schemes on unstructured polygonal grids by 
Kr\"oner and Rokyta, \cite{KroenerRokyta}, Kr\"oner, Noelle and Rokyta
\cite{KroenerNoelleRokyta}, Noelle \cite{Noelle1995}, and others.
Rohde \cite{Rohde} extended the method to weakly coupled systems
of conservation laws.

\subsection{Kinetic formulation and velocity averaging}
The kinetic formulation was introduced by Lions, Perthame,
and Tadmor \cite{LPT}. They show that there is a one-to-one
correspondence between weak entropy solutions of the scalar
conservation law (\ref{1_100}) and solutions of a linear
transport equation with source term, for which a certain nonlinear 
constraint holds. More precisely, one considers a ``density-like''
function $\rho$ depending on $(t,x)\in[0,t]\times\R^d$,
and on an additional variable $v\in\R$, 
which is a solution of the transport equation
\begin{equation}
\label{21_100}
\frac{\partial\rho}{\partial t}+f'(v)\cdot\nabla\rho
=\frac{\partial m}{\partial v}
\quad\mbox{in distributional sense}.
\end{equation}
Here $m$ is a nonnegative bounded measure. Equation
(\ref{21_100}) is supplemented with an assumption on the
structure of $\rho$. If $\chi\colon\R\rightarrow\R$ is defined
by
\begin{equation}
\label{21_200}
\chi(v|\alpha):=\left\{
\begin{array}{rcl}
+1 && \mbox{if $0<v<\alpha$} \\
-1 && \mbox{if $\alpha<v<0$} \\
 0 && \mbox{otherwise}
\end{array}
\right.
\end{equation}
for $\alpha\in\R$, then $\rho$ should have the form
\begin{equation}
\label{21_300}
\rho(t,x,v)=\chi(v|u(t,x))
\end{equation}
for some scalar function $u$. If $\rho$ is a
solution of the kinetic equation (\ref{21_100}) and satisfies
the nonlinear constraint (\ref{21_300}), then the function $u$
is the unique weak entropy solution of (\ref{1_100}).
Vice versa, if $u$ is the unique entropy solution of (\ref{1_100}),
then there exists a nonnegative bounded measure $m$ such that
the function $\rho$ defined by (\ref{21_300}) is a solution of
the kinetic equation (\ref{21_100}).

One big advantage of the kinetic formulation for scalar
conservation laws is the possibility to apply velocity averaging
lemmas to obtain regularity and compactness for (sequences of approximate) 
weak entropy solutions of (\ref{1_100}). One can show that the macroscopic
quantity $u$ has more regularity than the $\rho$ whose $v$-average
it is, see Golse, Lions, Perthame and Sentis \cite{Golse} and
DiPerna, Lions and Meyer \cite{DiPernaLionsMeyer}.

In Westdickenberg and Noelle \cite{West} we used the kinetic 
formulation together with the velocity averaging technique to
prove the convergence of a class of finite volume schemes for
scalar conservation laws in several space dimensions. We will
discuss this result in more detail in Section \ref{SS33} below.

We remark that already in 1984, Brenier \cite{Brenier}
constructed and analyzed an approximate evolution operator
for (\ref{1_100}) based on a kinetic decomposition of $u$.
In 1998 Perthame \cite{Perthame1998} used the kinetic formulation 
to study again uniqueness of weak entropy solutions and error 
estimates in the sense of Kuznetsov's approximation theory.
It is also possible to give a kinetic formulation for certain
systems of conservation laws, e.g. the isentropic Euler 
equations in 1-d. In Lions, Perthame, Tadmor \cite{LPT2} and 
Lions, Perthame, Souganidis \cite{LPS} 
this kinetic formulation was used together with 
compensated compactness to prove existence of weak entropy
solutions, thereby extending DiPerna \cite{DiPerna1983} and Chen's 
\cite{Chen} results mentioned above.  

\section{Some convergence results for finite volume schemes}
\label{S3}
In this section we give an overview over some of our own results
on the convergence of  finite volume schemes. These schemes 
define approximate solutions which are piecewise polynomial
on a given polygonal grid of $\R^d$, and in general discontinuous
at the cell interfaces. Integrating the conservation law (\ref{1_100})
over a cell $K$ with faces $e \subset \partial K$ from time
$t^n$ to time $t^{n+1}=t^n + \Delta t$ gives the
following update formula for the cell average $u_K^{n+1}$:
\begin{equation}
\label{3_100}
u_K^{n+1} = u_K^{n} - \frac{\Delta t}{|K|} \sum\limits_{e \subset \partial K}
|e| \, g_{K,e}^n.
\end{equation}
Here $g_{K,e}^n$ is a numerical flux which is consistent with 
$f(u)\cdot n_{K,e}$, where $n_{K,e}$ is the outward unit normal
of face $e$. This flux has to be conservative, satisfy some
upwinding property (e.g. an $E$-flux),
and it may be higher order accurate in space and time,
see e.g. \cite{KroenerNoelleRokyta,Noelle1995}. We consider
explicit time discretizations, which means that the fluxes $g_{K,e}^n$
can be computed directly from the cell averages $u_{K'}^n$ in
a neighborhood of $K$. Such discretizations are only stable
for timesteps $\Delta t$ which are so small that waves
originating from face $e=K\cap K'$, where the approximation is
discontinuous, do not cross the neighboring cells $K$ and $K'$
completely during one timestep. If one imposes this CFL-condition
(named after Courant, Friedrichs and Lewy \cite{CFL}),
then it is possible to prove  a discrete maximum principle.

\subsection{Discrete entropy inequalities}
\label{SS31}
The discrete maximum principle is enough to assure that
a subsequence converges weak$^*$ to some limit function,
but it does not guarantee that this limit is a
weak solution of the conservation law. The crucial
additional estimate needed to prove convergence
is a discrete entropy inequality of the form
\begin{equation}
\label{31_100}
\eta(u_K^{n+1}) - \eta(u_K^{n}) 
+ \frac{\Delta t}{|K|} \sum\limits_{e \subset \partial K}
|e| \, G_{K,e}^n \leq C h^{2\alpha}.
\end{equation}
Here $h$ is the maximal diameter of a cell of the computational grid,
$\alpha \in (\frac12,1]$ is a constant depending on the triangulation
and the polynomial reconstruction, the constant $C$ depends on
the entropy $\eta$ and its derivatives, and $G_{K,e}^n$ is a
numerical entropy flux which is consistent
with the entropy flux $q \cdot n_{K,e}$ and is closely related 
to the numerical flux $g_{K,e}^n$.

In \cite{KroenerNoelleRokyta} we proved an entropy inequality
of this type for a class of higher order schemes based on
the Lax-Friedrichs and the Engquist-Osher schemes,
and in \cite{Noelle1995} we could treat schemes which extend 
Godunov's exact Riemann solver to higher order accuracy
(see \cite{KroenerNoelleRokyta,Noelle1995} for the
definitions and the original references).

There are several points in the derivation of estimate 
(\ref{31_100}) which are worth recalling.
In \cite{Noelle1995} we first rewrite the update (\ref{3_100})
as a convex combination of one-dimensional contributions.
This is natural, since the numerical fluxes $g_{K,e}^n$
approximate one-dimensional fluxes $f \cdot n_{K,e}$
propagating in the normal direction to the cell-faces.
Adapting the weights of the decomposition
to the local wave speeds, our analysis admits larger time
steps than previous work, even for the classical first
order schemes in 1-d \cite{Tadmor1984}. 

Next we use Tadmor's decomposition \cite{Tadmor1984} of any $E$-flux
as a convex combination of Lax-Friedrichs and Godunov's fluxes.
If one uses a piecewise constant reconstruction, corresponding
to a first order accurate discretization, then the discrete
entropy inequality with $C = 0$ follows naturally from the
upwinding properties and the CFL-condition.
For higher order polynomial reconstructions,
we have to generalize Harten's TVD conditions discussed in
Section 2.1 above to multidimensional reconstructions.
The estimate for the higher order Lax-Friedrichs numerical flux follows
readily, but for the higher order Godunov's flux, which yields the best
resolution of shock discontinuities, the structure of the Riemann
solution has to be exploited in full detail, see \cite{Noelle1995}.
Once the one-dimensional discrete entropy inequalities
are obtained, convexity arguments yield the multidimensional
inequality (\ref{31_100}).


\subsection{Convergence and error estimates via $L^1$-contraction}
\label{SS32}
Equipped with the discrete maximum principle and entropy inequality
(\ref{31_100}), we can now prove convergence and error estimates.
In \cite{KroenerNoelleRokyta} we gave the first prove of convergence
for higher order accurate finite volume schemes on unstructured grids
for initial data which are merely in $L^\infty$, using DiPerna's
theory of measure-valued solutions. Subsequently, Kr\"oner and collaborators
designed practical algorithms along the lines of our convergence results.

In closely related work, Cockburn, Coquel and LeFloch \cite{CockburnCoquelLeFloch}
used Kuz\-net\-sov's approximation theory to prove convergence and error
estimates for BV initial data in 1994. In \cite{Noelle1996}, we could generalize
these error estimates, and simplify their proof,
using our discrete entropy inequalities. 

Even though inequality (\ref{31_100}) is the crucial estimate, 
careful work is required to finish the convergence proof.
In this process, many papers (including \cite{KroenerNoelleRokyta})
assume that the triangulation is regular, meaning that the ratio
of the outer diameter to the inner diameter of a cell is bounded above.
In \cite{Noelle1995}, we removed that restriction and treated
grids where the cells may become flat in the limit as
$h \to 0$. In \cite{Noelle1996}, we quantified the influence
of irregular grids on the convergence rate. An example in \cite{Noelle1995}
shows that the rate with which the grid may degenerate is optimal.


\subsection{Convergence via kinetic formulation}
\label{SS33}
In Westdickenberg and Noelle \cite{West} we gave a new convergence
proof for sequences of approximate solutions of the scalar conservation
law, produced by higher order finite volume schemes on unstructured 
grids in several space dimensions. We did not use
the theory of measure-valued solutions, but relied on the kinetic
formulation and velocity averaging lemmas instead. The key tools
in our analysis were once more the discrete entropy condition (\ref{31_100})
and the following compactness theorem.
\begin{theorem}
\label{T21_100}
Let $1<p\leq 2$ and $0<\gamma<1$. Fix a compact set $K\subset\R$
and assume that there are sequences $\{\rho^\epsilon\}_\epsilon, 
\{m^\epsilon\} _\epsilon$, and $\{\pi^\epsilon\}_\epsilon$ uniformly 
bounded in $L^p(\R^{d+1}\times K), L^1(K,\mathcal{M}(\R^{d+1}))$,
and $L^1(K, C_\gamma^*(\R^{d+1}))$, resp., such that
\begin{equation}
\label{21_400}
\frac{\partial \rho^\epsilon}{\partial t}+f'(v)\cdot\nabla
\rho^\epsilon=\frac{\partial m^\epsilon}{\partial v}+\pi^\epsilon
\quad\mbox{in distributional sense}
\end{equation}
($(t,x)\in\R^{d+1}$, $v\in K$). 
If the following nondegeneracy condition holds
\begin{equation}
\label{21_500}
\sup_{(\tau,\xi)\in\R^{d+1}}\mathrm{meas}\bigl\{v\in\Lambda\colon
   \tau+f'(v)\cdot\xi=0\bigr\}=0,
\end{equation}
then the sequence $\{z^\epsilon\}_\epsilon$ belongs to a compact
subset of $L^1_{loc}(\R^{d+1})$, where
\begin{equation}
\label{21_600}
z^\epsilon:=\int_K\rho^\epsilon(\cdot,v)\,dv.
\end{equation}
\end{theorem}
Here, $\mathcal{M}(\R^{d+1})$ is the space of bounded Radon measures,
and $C_\gamma^*(\R^{d+1})$ is the topological dual of the space of
H\"older continuous functions.
This theorem is a variant of the compactness theorem in \cite{LPT}.
Assumption (\ref{21_500}) means that the flux function must be 
nonlinear enough. The regularizing effect stated in Theorem 
\ref{T21_100} does not exist for linearly degenerate problems
(e.g. for advection equations).
Note that the $m^\epsilon$ on the right hand side of (\ref{21_400})
is a bounded measure,
as we should expect from the kinetic formulation
(\ref{21_100}). The quantity $\pi^\epsilon$ is an error term
(e.g. the numerical residual). It is measured in a function space 
with negative regularity, hence in a rather weak topology.

To prove Theorem \ref{T21_100} we decompose $z^\epsilon$ into 
two parts and show that
one part can be made arbitrarily small in $L^1_{loc}$ using the 
nondegeneracy of the flux (\ref{21_500}), while the other part has
some positive regularity in terms of Sobolev (in fact Besov) spaces
and is therefore strongly compact in $L^1_{loc}$.
Applying Theorem \ref{T21_100} to numerical approximations
we obtain the strong compactness needed to pass to
the limits in nonlinear quantities as explained in Section \ref{Intro}. 

\section{Related work and open problems}

In \cite{KroenerOhlberger2000,Ohlberger2001} Kr\"oner and Ohlberger
derived a posteriori error estimates based on Kuznetsov's
approximation theory and recent related work of Gallouet, Herbin,
Chainais-Hillairet et al. \cite{Chainais1999}.
Ohlberger \cite{Ohlberger2001} developed a fully adaptive,
implicit finite volume scheme
for scalar convection-reaction-diffusion equations based
on these estimates.

In a recent preprint, Hwang and Tzavaras \cite{HwangTzavaras} used the kinetic 
formulation together with the velocity averaging technique to study 
the convergence of approximate solutions of scalar conservation laws
that are produced by either a relaxation method or by a 
diffusion-dispersion approximation.

To conclude this note, we would like to sketch some 
important open problems concerning the convergence of
numerical approximations.
The convergence rate of $h^{1/2}$ proved by Kuznetsov for
the Lax-Friedrichs scheme is only optimal
for first order approximations of discontinuous
solutions of linear advection equations. The
rate obtained so far for unstructured grids
is only $h^{1/4}$. A first attempt to overcome this
barrier may be found in \cite{CockburnGremaudYang1998}.
For strictly convex scalar conservation laws one
expects a rate of $h \log h$, but this conjecture is so far
only backed up by numerical experiments.
For schemes which are formally higher order accurate
one expects higher convergence rates away from 
discontinuities. For strictly convex scalar conservation laws
a result in this direction is proven in \cite{EngquistSjoegreen1998}.
It is also shown there that for systems of conservation laws,
characteristics of one family crossing a numerical shock layer
belonging to a different family may
carry first order errors into the smooth postshock
region, so the convergence rate deteriorates there even for
higher order schemes. In ongoing work G.~Kreiss et al.
are beginning to design special numerical schemes which
attempt to filter out this pollution effect.


\begin{thebibliography}{8.}
\bibitem{BouchutPerthame}
      {\sc F.~Bouchut and B.~Perthame},
      {\em Kruzkov's estimates for scalar conservation laws revisited},
      Trans. AMS~350 (1998), pp.~2847--2870.
\bibitem{Brenier}
      {\sc Y.~Brenier},
      {\em Averaged multivalued solutions for scalar conservation laws.}
      SIAM J. Numer. Anal. 21 (1984), pp.~1013--1037.
\bibitem{Bressan}
      {\sc A.~Bressan},
      {\em Hyperbolic systems of conservation laws. The
        one-dimensional Cauchy problem},
      Oxford Lecture Series in Mathematics and its Applications 20,
      Oxford University Press, Oxford, 2000.
\bibitem{Chainais1999}
      {\sc C.~Chainais-Hillairet},
      {\em Finite volume schemes for a nonlinear hyperbolic equation.
	Convergence towards the entropy solution and error estimate},
      M2AN Math. Model. Numer. Anal.~33 (1999), pp.~129--156.
\bibitem{Chen}
      {\sc G.-Q.~Chen},
      {\em Convergence of the Lax-Friedrichs scheme for isentropic 
        gas dynamics III},
      Acta Math. Sci.~6 (1986), pp.~75--120.	
\bibitem{ChenDuTadmor}
      {\sc G.-Q.~Chen, Q.~Du and E.~Tadmor},
      {\em Spectral viscosity approximations to multidimensional
        scalar conservation laws},
      Math. Comp. 61 (1993), pp.~629--643.
\bibitem{CockburnCoquelLeFloch}
      {\sc B.~Cockburn, F.~Coquel and P.~LeFloch},
      {\em An error estimate for finite volume methods for 
        multidimensional conservation laws},
      Math. Comp. 63 (1994), pp.~77--103.
\bibitem{CockburnGremaudYang1998}
      {\sc B.~Cockburn, P.-A.~Gremaud and J.X.~Yang},
      {\em A priori error estimates for numerical methods for scalar
        conservation laws. III. Multidimensional flux-splitting
        monotone schemes on non-Cartesian grids.},
      SIAM J.~Numer. Anal. 35 (1998), pp.~1775--1803.
\bibitem{ConwaySmoller}
      {\sc E.~Conway and J.~Smoller},
      {\em Global solutions of the Cauchy problem for quasilinear 
        first-order equations in several space variables},
      Comm. Pure Appl. Math.~19 (1966), pp.~95--105.
\bibitem{CoquelLeFloch}
      {\sc F.~Coquel and P.~LeFloch},
      {\em Convergence of finite difference schemes for conservation
        laws in several space dimensions: the corrected antidiffusive
        flux approach},
      Math.~Comp.~57 (1991), pp.~169--210.
\bibitem{CFL}
      {\sc R.~Courant, K.O.~Friedrichs and H.~Lewy},
      {\em \"Uber die partiellen Differentialgleichungen der
         mathematischen Physik},
      Math. Annal.~100 (1928), pp.~32--74.
\bibitem{Dafermos}
      {\sc C. Dafermos},
      {\em Hyperbolic conservation laws in continuum physics},
      Grund\-lehren der Mathematischen Wissenschaften
        (Fundamental Principles of Mathematical Sciences) 325,
        Springer, Berlin, 2000.
\bibitem{DiPerna1983}
      {\sc R.DiPerna},
      {\em Convergence of the viscosity method for isentropic 
        gas dynamics},
      Comm. Math. Phys.~91 (1983), pp.~1--30.
\bibitem{DiPerna1985}
      {\sc R.~DiPerna},
      {\em Measure-valued solutions to conservation laws},
      Arch. Rat. Mech. Anal. 88 (1985), pp.~223-270.
\bibitem{DiPernaLionsMeyer}
      {\sc R.J.~DiPerna, P.-L.~Lions and Y.~Meyer},
      {\em $L^p$ regularity of velocity averages},
      Ann. I.H.P. Analyse Nonlin\'eaire, 8 (1991), pp.~271--287.
\bibitem{EngquistSjoegreen1998}
      {\sc B.~Engquist and B.~Sj\"ogreen},
      {\em The convergence rate of finite difference schemes in the
	presence of shocks},
      SIAM J. Numer. Anal.~35 (1998), pp.~2464--2485.
\bibitem{Glimm}
      {\sc J.~Glimm},
      {\em Solutions in the large for nonlinear hyperbolic
        systems of equations},
      Comm. Pure Appl. Math.~18 (1965), pp.~697--715.
\bibitem{Golse}
      {\sc F.~Golse, P.-L.~Lions, B.~Perthame and R.~Sentis},
      {\em Regularity of the moments of the solution of a
        transport equation},
      J.~Funct. Anal. 76 (1988), pp.~110--125.
\bibitem{Harten}
      {\sc A.~Harten},
      {\em High resolution schemes for hyperbolic conservation laws},
      J.~Comp. Phys.~49 (1983), pp.~357--393.
\bibitem{Hopf}
      {\sc E.~Hopf},
      {\em The partial differential equation 
        $u_t+uu_x=\mu u_{xx}$},
      Comm. Pure Appl. Math.~3 (1950), pp.~201--230.
\bibitem{HwangTzavaras}
      {\sc S.~Hwang and A.E.~Tzavaras},
      {\em Kinetic decomposition of approximate solutions to 
        conservation laws: application to relaxation and 
	diffusion-dispersion approximations},
      Preprint, University of Wisconsin-Madison, 2002.	
\bibitem{KroenerNoelleRokyta}
      {\sc D.~Kr\"oner, S.~Noelle and M.~Rokyta},
      {\em Convergence of higher order upwind finite volume schemes
        on unstructured grids for scalar conservation laws in several
        space dimensions},
      Numer. Math. 71 (1995), pp.~527--560.
\bibitem{KroenerOhlberger2000}
      {\sc D.~Kr\"oner and M.~Ohlberger},
      {\em A posteriori error estimates for upwind finite volume schemes
        for nonlinear conservation laws in multidimensions},
      Math. Comp.~69 (2000), pp.~25--39.
\bibitem{KroenerRokyta}
      {\sc D.~Kr\"oner and M.~Rokyta},
      {\em Convergence of upwind finite volume schemes for scalar
        conservation laws in two dimensions},
      SIAM J.~Numer. Anal. 31 (1994), pp.~324--343.
\bibitem{Kruzkov}
      {\sc S.N.~Kruzkov},
      {\em First order quasilinear equations in several independent
        variables},
      Math. Sb. 123 (1970), pp.~228--255; 
      English transl. in Math. USSR Sbornik~10 (1970), pp.~217--243.
\bibitem{Kuznetsov}
      {\sc N.N.~Kuznetsov},
      {\em Accuracy of some approximate methods for computing the weak
        solutions of a first-order quasi-linear equation},
      USSR Comp. Math. Math. Phys. 16/6 (1976), pp.~105--119.
\bibitem{Ladyzhenskaya}
      {\sc O.A.~Ladyzhenskaya},
      {\em On the construction of discontinuous solutions
        of quasi-linear hyperbolic equations as limits of
	solutions of the corresponding parabolic equations
	when the ``coefficient of viscosity'' tends towards
	zero} (in Russian),
      Trudy Moscov.~Mat.~Obsc. 6 (1957), pp. 465--480.
\bibitem{Lax1954}
      {\sc P.~Lax},
      {\em Weak solutions of nonlinear hyperbolic equations and
        their numerical approximation},
      Comm. Pure Appl. Math.~7 (1954), pp.~159--193.
\bibitem{Lax1973}
      {\sc P.~Lax},
      {\em Hyperbolic Systems of Conservation Laws and the Mathematical
        Theory of Shock Waves},
      SIAM, Philadelphia, 1973.
\bibitem{Lax1971}
      {\sc P.~Lax},
      {\em Shock waves and entropy},
      in {\em Contributions to nonlinear functional analysis},
        E.A.~Zarantonello Ed.,
        Academic Press, New York, 1971.
\bibitem{LPS}
      {\sc P.-L.~Lions, B.~Perthame and P.E.~Souganidis},
      {\em Existence and stability of entropy solutions for the
        hyperbolic systems of isentropic gas dynamics in Eulerian
	and Lagrangian coordinates},
      Comm. Pure Appl. Math.~49 (1996), pp.~599-638.	
\bibitem{LPT}
      {\sc P.-L.~Lions, B.~Perthame and E.~Tadmor},
      {\em A kinetic formulation of multidimensional scalar
        conservation laws and related equations},
      J.~AMS, 7 (1994), pp.~169--191.
\bibitem{LPT2}
      {\sc P.-L.~Lions, B.~Perthame and E.~Tadmor},
      {\em Kinetic formulation of the isentropic gas dynamics 
        and $p$-systems},
      Comm. Math. Phys.~163 (1994), pp.~415--431.
\bibitem{Liu}
      {\sc T.-P.~Liu},
      {\em The entropy condition and the admissibility of shocks},
      J.~Math. Anal. Appl.~53 (1976), pp.~78--88.
\bibitem{Murat}
      {\sc F.~Murat},
      {\em Compacit\'e par compensation},
      Ann. Sc.~Norm. Sup. Pisa 5 (1978), pp.~489--509.
\bibitem{Noelle1995}
      {\sc S.~Noelle},
      {\em Convergence of higher order finite volume schemes on
        irregular grids},
      Adv. Comp. Math. 3 (1995), pp.~197--218.
\bibitem{Noelle1996}
      {\sc S.~Noelle},
      {\em A note on entropy inequalities and error estimates for
        higher order accurate finite volume schemes on irregular
        grids},
      Math. Comp. 65 (1996), pp.~1155--1163.
\bibitem{Ohlberger2001}
      {\sc M.~Ohlberger},
      {\em A posteriori error estimates for vertex centered finite volume
        approximations of convection-diffusion-reaction equations},
      M2AN Math. Model. Numer. Anal.~35 (2001), pp.~355--387.
\bibitem{Oleinik1957}
      {\sc O.A.~Oleinik},
      {\em Discontinuous solutions of nonlinear differential equations},
      Usp. Mat. Nauk.~12 (1957), pp.~3--73; 
      English transl. in AMS Transl.~26 (1963), pp.~95--172.
\bibitem{Osher}
      {\sc S.~Osher},
      {\em Convergence of generalized MUSCL schemes},
      SIAM J.~Numer. Anal. 22 (1985), pp.~947--961.
\bibitem{OsherTadmor}
      {\sc S.~Osher and E.~Tadmor},
      {\em On the convergence of difference approximations
      to scalar conservation laws},
      Math. Comp. 50 (1988), pp.~19--51.
\bibitem{Perthame1998}
      {\sc B.~Perthame},
      {\em Uniqueness and error estimates in first order quasilinear 
        conservation laws via the kinetic entropy defect measure},
      J.~Math. Pures Appl.~77 (1998), pp.~1055--1064.
\bibitem{Riemann}
      {\sc B.~Riemann},
      {\em \"Uber die Fortpflanzung ebener Luftwellen von endlicher
        Schwin\-gungsweite},
      Abh. d. K\"onigl. Ges. d. Wiss. z. G\"ottingen, Bd.~8 (1858/59) (Math. Cl.),
      pp.~43--65.
\bibitem{Rohde}
      {\sc C.~Rohde},
      {\em Upwind finite volume schemes for weakly coupled
        hyperbolic systems of conservation laws in 2D},
      Numer. Math. 81 (1998), pp.~85--123.
\bibitem{Sanders}
      {\sc R.~Sanders},
      {\em On convergence of monotone finite difference schemes with
        variable spatial differencing},
      Math. Comp. 40 (1983), pp.~91--106.
\bibitem{Smoller}
      {\sc J.~Smoller},
      {\em Shock Waves and Reaction-Diffusion Equations},
      Springer, New York, 1983.
\bibitem{Szepessy1989a}
      {\sc A.~Szepessy},
      {\em An existence result for scalar conservation laws using
        measure valued solutions},
      Comm. Part. Diff. Equations 14 (1989), pp.~1329--1350.
\bibitem{Szepessy1989b}
      {\sc A.~Szepessy},
      {\em Convergence of a shock-capturing streamline diffusion 
        finite element method for scalar conservation laws in two 
	space dimensions},
      Math. Comp. 50 (1989), pp.~527--545.
\bibitem{Tadmor1984}
      {\sc E.~Tadmor},
      {\em Numerical viscosity and the entropy condition
        for conservative difference schemes},
      Math. Comp.~43 (1984), pp.~369--381.
\bibitem{Tadmor1991}
      {\sc E.~Tadmor},
      {\em Local error estimates for discontinuous solutions of nonlinear
        hyperbolic equations},
      SIAM J.~Numer. Anal. 28 (1991), pp.~811--906.
\bibitem{Tartar}
      {\sc L.~Tartar},
      {\em Compensated compactness and applications to partial differential
        equations},
      in {\em Research Notes in Mathematics 39}, Nonlinear Analysis and
        Mechanics, Herriot-Watt Symposium, R.J.~Knopps Eds.,
        Pittman Press, 1975.
\bibitem{Vila}
      {\sc J.P.~Vila},
      {\em Convergence and error estimates in finite volume schemes for
        general multi-dimensional scalar conservation laws I. Explicit
        monotone schemes},
      RAIRO Anal. Num\'er. 28 (1994), pp.~267--295.
\bibitem{Volpert}
      {\sc A.I.~Vol'pert},
      {\em The space BV and quasilinear equations},
      Math. USSR Sbornik~2 (1967), pp.~225-267.
\bibitem{Vvedenskaya1956}
      {\sc N.D.~Vvedenskaya},
      {\em Solution of the Cauchy-problem for a nonlinear equation
        with discontinuous initial conditions by the method of
        finite differences},
      Dokl. Akad. Nauk SSSR~111 (1956), pp.~517-520. (in Russian)
\bibitem{West}
      {\sc M.~Westdickenberg and S.~Noelle},
      {\em A new convergence proof for finite volume schemes
        using the kinetic formulation of conservation laws},
      SIAM J.~Num. Anal. 37 (2000), pp.~742--757.
\end{thebibliography}
\end{document}